\documentclass[twoside,leqno]{article}

\usepackage[letterpaper]{geometry}
\usepackage{siamproceedings}
\usepackage[T1]{fontenc}
\usepackage{amsfonts}
\usepackage{graphicx}
\usepackage{epstopdf}
\usepackage{enumitem}
\usepackage{algorithmic}
\usepackage{amssymb}
\usepackage{amsmath}

\ifpdf
  \DeclareGraphicsExtensions{.eps,.pdf,.png,.jpg}
\else
  \DeclareGraphicsExtensions{.eps}
\fi

\usepackage{comment}
\renewcommand{\phi}{\varphi}
\renewcommand{\leq}{\le}
\renewcommand{\geq}{\ge}
\newcommand{\one}{\mathbf{1}}

\newcommand{\eps}{\varepsilon}

\newcommand{\s}{\sigma}

\newcommand{\cE}{\mathcal E}

\newcommand{\cH}{\mathcal H}

\newcommand{\cI}{\mathcal I}
\newcommand{\1}{\mathbb{1}}

\newcommand{\Per}{\operatorname{Per}}
\def\1{\mathbbm{1}}

\def\g{{\gamma}}
\renewcommand{\le}{\leqslant}
\renewcommand{\ge}{\geqslant}
\renewcommand{\P}{\mathbb{P}}
\newcommand{\R}{\mathbb R}

\newcommand{\Z}{\mathbb Z}

\newcommand{\N}{\mathbb N}

\newcommand{\EE}{\mathbb E}

\newcommand{\PP}{\mathbb P}
\newcommand{\bS}{\mathbb S}
\newcommand{\la}{\langle}
\newcommand{\ra}{\rangle}
\newcommand{\Vol}{\mathrm{Vol}}

\newsiamremark{remark}{Remark}
\newsiamremark{hypothesis}{Hypothesis}

\crefname{hypothesis}{Hypothesis}{Hypotheses}
\newsiamthm{claim}{Claim}
\newsiamthm{conjecture}{Conjecture}

\usepackage{amsopn}

\begin{document}

\title{\Large Probabilistic combinatorics at exponentially small scales}
    \author{Julian Sahasrabudhe\thanks{University of Cambridge (\email{jdrs2@cam.ac.uk}, \url{https://www.dpmms.cam.ac.uk/~jdrs2/}).}
    }

\date{}

\maketitle


\fancyfoot[R]{\scriptsize{Copyright \textcopyright\ 20XX by SIAM\\
Unauthorized reproduction of this article is prohibited}}





\begin{abstract}
 In many applications of the probabilistic method, one looks to study phenomena that occur ``with high probability''. More recently however, in an attempt to understand some of the most fundamental problems in combinatorics, researchers have been diving deeper into these probability spaces, and understanding phenomena that occur at much smaller probability scales. Here I will survey a few of these ideas from the perspective of my own work in the area. 

\end{abstract}

\maketitle

\section{Introduction} 
In many applications of the probabilistic method one aims to show the existence of an object not by constructing the object directly, but rather by setting up a probability space and showing that the one can draw the object of interest with non-zero probability. While this might at first seem like a trivial recasting of the problem, it has proven to be incredibly powerful change of perspective which has come to dominate the field in recent years. 

Classically, in such applications of the probabilistic method, one studies events that occur not just with non-zero probability, but with probability close to one. A state of affairs which has been cheekily sloganized as ``the probabilistic method finds the hay in the haystack''. In this survey we will touch on a few topics where we are forced to go far beyond these typical behaviors and study phenomena that occur at exponentially small probability scales - right at the edge where probability is still useful. Perhaps surprisingly, there is quite a bit to be said in these cases. 

Rather than attempt to properly survey this topic, which I am not properly placed to do,  I will instead use this theme to tie together a few strands of my own recent work and to highlight some of my thinking on the topics where I have had some success. 

\vspace{2mm}

\emph{Flat littlewood polynomials.} We begin this survey in Section~\ref{sec:discrepancy} where we introduce the area of  \emph{discrepancy theory} before sketching how some of these ideas came to inform the thinking behind a recent result of myself, joint with Paul Balister, B\'{e}la Bollob\'{a}s, Robert Morris and Marius Tiba, in the area of harmonic analysis \cite{balister2020}. Here we used ideas from discrepancy theory to construct so called \emph{flat Littlewood polynomials}, thus resolving an old conjecture of Littlewood.

A Littlewood polynomial is a polynomial with all coefficients in $\{-1,1\}$. Inspired by a question of Erd\H{o}s \cite{E57}, Littlewood went on to consider the following question regarding about how ``flat'' the modulus of such polynomials can be.  If $P$ is a degree $n$ Littlewood polynomial then, by a simple application of Parseval, we see that we must have that $\max_{z : |z| =1 } |P(z)| \geq \sqrt{n+1}$. Littlewood conjectured that there exist such polynomials that are ``flat'' as possible, in the sense that $|P(z)| = \Theta(\sqrt{n})$ for all $z$ with $|z|=1$. 

Using tools from probability and discrepancy theory, we solved this conjecture by showing that there are constants\/ $c_1,c_2 >  0$ so that, for all\/ $n \ge 2$,
 there exists a Littlewood polynomial\/ $P(z)$ of degree $n$ with 
\begin{equation}
c_1 \sqrt{n} \le |P(z)| \le c_2 \sqrt{n}
 \end{equation}
 for all\/ $z \in \mathbb{C}$ with\/ $|z| = 1$. This set us up to discuss some beautiful open problems in discrepancy theory and some exciting recent developments.

\vspace{2mm}

\emph{Sphere packing in high dimensions.} In Section~\ref{sec:processes} we will go on to discuss a technique which can used to ``access'' these events at very small probability scales; via a random process that is biased towards the event of interest. We will use this thinking to frame the recent work of myself, joint with Marcelo Campos, Matthew Jenssen and Marcus Michelen, where we construct sphere packings and spherical codes in large dimension. In particular, we show the following new lower bounds on the classical sphere packing problem, giving the first asymptotically growing improvement since the work 1947 work of Rogers \cite{rogers1947existence}.

Let $\theta(d)$ be the maximum proportion of $\mathbb{R}^d$ which is covered by non-overlapping identical spheres. As the dimension $d$ tends to infinity, we have 
    \begin{equation}\label{eq:our-spheres}
    \theta(d) \geq (1-o(1)) \frac{d\log d}{2^{d+1}}\, .
    \end{equation}
    One interesting aspect of our sphere packings is that they are both ``dense'' and highly ``disordered''. Such sphere packings arise in the physics literature and were conjectured to exist by Parisi and Zamponi \cite{parisi2006amorphous}. It moreover seems reasonable to conjecture that the bound~\eqref{eq:our-spheres} is actually sharp, up to constant factors, for such ``highly disordered'' sphere packings. 
    
This method also naturally adapts to improve the lower bounds on spherical codes in large dimension. Let $A(d,\theta)$ be the maximum proportion of the sphere $\bS^{d-1}$ that can be covered by non-overlapping spherical caps of radius $\theta$. If $\theta \in (0,\pi/2),$ then
\[ A(d,\theta)\geq (1-o(1)) \frac{d\log d}{2 s_d(\theta)},\]
where $s_d(\theta)$ denotes the normalized volume of a spherical cap. 

We will then go on to discuss a further beautiful advance of Klartag \cite{klartag2025} who has construed even denser sphere packings in large dimension by studying properties of random lattices at exponentially small scales. He proved that
    \[
    \theta(d) \geq (c-o(1)) d^2/2^{d}\, ,    \] 
for some $c>0$. Interestingly, it seems reasonable to conjecture that Klartag's bound is also sharp, up to constant factors, in the case of \emph{lattice} packings.

Despite these advances, there still remain exponential gaps between the upper and lower bounds on the quantities $\theta(d)$ and $A(d,\theta)$ as $d\rightarrow \infty$ and it seems to be an incredibly inciting question to improve either by an exponential factor. 

\vspace{2mm}

\emph{Random matrix theory at exponentially small scales. } In the final section we shall keep with our theme of high dimensional geometry and move on to discuss phenomena that occur at exponentially small scales in the setting of random matrix theory. This discussion will center around a recent result of the author in joint work with Campos, Jenssen and Michelen \cite{RSM2} on the singularity probability of a random symmetric matrix.
We showed that if $A_n$ is drawn uniformly at random from the $n\times n$ symmetric matrices with entries in $\{-1,1\}$ then 
\begin{equation}\label{eq:Mn-singular} \PP\big( \det(A_n) = 0 \big) \leq e^{-cn} ,\end{equation}
where $c>0$ is an absolute constant.

As it happens, this will open us to discuss a handful of other results on related topics - the study of the distribution of the least singular value and  estimates on the repulsion and clustering of the eigenvalues of random symmetric matrices. We will, in particular, highlight several of the techniques developed in service of \eqref{eq:Mn-singular}. 

\section{Beck, Spencer and Flat Littlewood polynomials}\label{sec:discrepancy}

Motivated by generalizing his famous theorem on $3$-term progressions in dense sets of integers \cite{roth1953}, Roth considered the following natural problem on two-colourings of an interval. Let $n$ be a large integer and let $f$ be an arbitrary function\footnote{Here and throughout we use the notation $[n] = \{ 1,\ldots, n\}$.}  
$ f : [n] \rightarrow \{-1,1\}$, which we extend to all of $\mathbb{N}$ by setting $f(x) = 0 $ for all $x>n$. Define the \emph{discrepancy} of $f$, denoted $D(f)$, to be the maximum sum in absolute value along an arithmetic progression;
\[ D(f) = \max_{a,b} \bigg| \sum_{x} f(ax+b) \bigg|.\]
Roth proved \cite{roth1964} that $D(f) \geq cn^{1/4}$ for all such functions $f$
and conjectured that this could be improved to $n^{1/2-o(1)}$. This was then disproved by Sark\H{o}zy who's result was improved a few more times before Beck \cite{beck1981roth} introduced his famous partial coloring method to resolve the question up to logarithms, showing that Roth's original bound was just about sharp
\[ \min_f D(f) = n^{1/4+o(1)}.\]
Some years later, the general utility of this beautiful idea was noticed by Spencer \cite{spencer1985} who used Beck's partial colouring idea to prove Spencer's famous theorem in discrepancy theory.

\begin{theorem}\label{thm:spencer}
 Let $A$ be a $n\times n$ matrix with all entires bounded in absolute value by one. Then there exists a vector $x \in \{-1,1\}^n$ for which 
 \begin{equation}\label{eq:spencer} \|Ax\|_{\infty} \leq C \sqrt{n},\end{equation}
 for some absolute constant $C>0$.
\end{theorem}

What is remarkable about this result is that if one takes a suitably typical matrix\footnote{Note we have to have some typicality condition here. For example if $A$ has rank $1$, then \eqref{eq:spencer} will hold with constant probability.} and a \emph{random} $x \in \{-1,1\}^n$ one will have 
\[\|Ax\|_2 = \Theta\big(\sqrt{n\log n}\big),\]
with high probability. Moreover, one expects that it is exponentially unlikely that a random vector $x \in \{-1,1\}^n$ satisfies 
$\|Ax\|_{\infty} \leq C\sqrt{n}$. This is because for each $i$ we have $|(Ax)_i| \leq C\sqrt{n}$ with constant probability and one expects that each entry is roughly independent (again assuming $A$ is appropriately non-degenerate). Thus to find the  vector $x \in \{-1,1\}^n$ guaranteed by Theorem~\ref{thm:spencer}, ``standard'' probabilistic arguments are not available. Instead one can adapt the simple but ingenious differencing method of Beck to this setting. 

For this, we iterate the following lemma that finds a vector $x \in \{-1,0,1\}^n$ with not to many zero entries. 

\begin{lemma}\label{lem:recolouring}
Let $A$ be a $n\times n$ matrix with complex entries, bounded in absolute value by one. Then there exists a vector $x \in \{-1,0,1\}^n$ with at least $n/4$ non-zero terms for which 
 \begin{equation} \|Ax\|_{\infty} \leq C \sqrt{n},\end{equation}
 for some absolute constant $C>0$.
\end{lemma}
\begin{proof}[Proof sketch]
To prove this consider all images $\{ Ay : y \in \{0,1\}^n\}$. One then uses a (careful) application of the pigeon-hole principle to show that there is a subset $S \subset \{0,1\}^n$ with $|S| \geq 2^{(1-\eps) n}$, for which  $\|y-y'\|_{\infty}\leq C\sqrt{n}$ for all $y,y' \in S$. Since each vector $y \in \{0,1\}$ has at most $\sum_i^{d} \binom{n}{i}$ vectors with hamming distance $\leq d$, there must be a pair $y,y'\in S$ with the property that $x= y-y'$ has at least $d = n/4$ non-zero entries. 
\end{proof}

\subsection{Flat Littlewood polynomials}
Interestingly these fundamental questions of discrepancy theory connect with an old and famous conjecture of Littlewood on ``Littlewood polynomials''. Indeed, say that a polynomial $P(z)$ of degree $n$, is a \emph{Littlewood polynomial} if 
\[ P(z) = \sum_{k=0}^n \eps_k z^k \qquad \text{ where } \qquad \eps_k \in \{-1,1\} . \] 

The study of such polynomials has a long and fascinating history (see~\cite{B02} or~\cite{M17}), with their study going back to the work of Hardy and Littlewood~\cite{HL16} over 100 years ago, but was also popularized by the the work of Bloch and P\'olya~\cite{BP32}, the work of Littlewood and Offord~\cite{LO38, LO48}, and others~\cite{EO56,SZ54}. 


While there are many beautiful questions about such polynomials, we are interested here in a problem asked by Erd\H{o}s~\cite{E57} in 1957 which was then taken up and extensively studied by Littlewood~\cite{L61,L62a,L66a,L66} in a series of papers on the extremal properties of polynomials with restricted coefficients. For this we say a (sequence) of Littlewood polynomials $P_n(z)$ ($P_n$ of degree $n$) is \emph{flat} if there exist constants $c_1,c_2 >0$ such that\footnote{As discussed above, we recall that $\sqrt{n}$ is the natural normalization since 
$(2\pi)^{-1}\int |P_n(z)|^2 = \sum_k |\eps_k|^2 = n+1$, by Parseval. }
\begin{equation}\label{eq:def-flat-polys} c_1\sqrt{n} \leq |P_n(z)| \leq c_2\sqrt{n}. \end{equation}
In particular, Littlewood conjectured in~\cite{L66} that flat Littlewood polynomials exist. 

To see a first connection with discrepancy, we note that if one relaxes \eqref{eq:def-flat-polys} to the one-sided bound,
\begin{equation}\label{eq:rudin-shapiro} \max_{z : |z| = 1} |P(z)| \leq C\sqrt{n}. \end{equation}
we arrive at the natural ``analogue'' of \eqref{eq:spencer} in the setting of polynomials. Actually, the one-sided problem was solved in the 1950s by Rudin \cite{rudin1959} and Shapiro \cite{shapiro1951} who introduced the famous ``Rudin-Shapiro'' polynomials. Interestingly, one of the first applications that Spencer gave of his new theorem was to give a very different proof of the existence of polynomials satisfying \eqref{eq:rudin-shapiro}. To see the connection, we briefly sketch the idea.

\begin{theorem} For every $n \geq 2$ there is a Littlewood polynomial $P_n$ with $\deg(P_n) = n$ which satisfies \eqref{eq:rudin-shapiro}, for some absolute $C>0$.
\end{theorem}
\begin{proof}[Proof sketch] For $\theta \in [0,2\pi]$, define the vector
\[ v(\theta) = \big( 1,e^{i\theta},e^{2i\theta},\ldots, e^{ni\theta}\big). \]
Note if we apply Spencer's theorem to $v(2\pi j/n)$, for $j = 1,\ldots,n$, we can ensure the polynomial is small at the $n$th roots of unity. But this is not quite enough. To get around this we can also use Spencers theorem to ensure that all of the derivatives $P^{(j)}$ are also well controlled at roots of unity. While this might seem like we are putting infinitely many constraints on the sequence $\{-1,1\}^n $ it is actually okay since we need weaker and weaker control on each derivative. Then, via Taylor's theorem we can ensure that the polynomial is well behaved everywhere. \end{proof}

To attack Littlewood's conjecture, one might hope for a strengthening of Theorem~\ref{thm:spencer} that also provides a \emph{lower bound} on each entry. However it is not hard to see that this is not possible. Consider the $n+1$ vectors $a^{(0)},\ldots,a^{(n)}$ where we define 
\[ a^{(i)} = (-1,\ldots,-1,1,\ldots 1)\] to have exactly $i$ $1$s and and $n-i$ entries equal to $-1$. So given any $x \in \{-1,1\}^n$,  we can assume $\la x, a^{(0)} \ra \geq 0$ and thus $\la x, a^{(n)} \ra  = -\la x, a^{(1)} \ra\leq 0$. Now note that for each $i< n$, 
\[\big|\la a^{(i+1)},x \ra - \la a^{(i)},x \ra\big| \leq 2\]
and therefore there is always a vector $a^{(i)}$ for which $|\la x,a^{(i)}\ra| \leq 2$. Thus we see that any matrix $A$ which is defined by any $n$ vectors among the $\{ a^{(i)} \}_i$ provides a counterexample to this potential strengthening of the Spencer's theorem.

In joint work with Bollob\'{a}s, Balister, Morris and Tiba, we showed that we could in fact get around this obstacle in the context of flat Littlewood polynomials and use discrepancy theory methods to prove the existence of these polynomials.

\begin{theorem}\label{thm:flat-littlewood}
 For all\/ $n \ge 2$, there exists a Littlewood polynomial\/ $P(z)$ of degree $n$ with 
 \begin{equation}
  c_1 \sqrt{n} \le |P(z)| \le c_2 \sqrt{n}
 \end{equation}
 for all\/ $z \in \mathbb{C}$ with\/ $|z| = 1$. Here $c_1,c_2>0$ are absolute constants. 
\end{theorem}

In what follows we sketch a proof of Theorem~\ref{thm:flat-littlewood}.

\subsection{Sketch of the construction}
 To prove Theorem~\ref{thm:flat-littlewood} we may assume that the degree is $4n$. We then write $z = e^{ix}$ and multiply by a phase $e^{-2ixk}$ so that we obtain a function 
 \[ f(x) =   \sum_{k = -2n}^{2n} \eps_k e^{ix k } . \]  Thus our goal is to choose the $\eps_k \in \{-1,1\}$ so that $f(x) = \Theta( n^{1/2})$, for all $x$. We now constrain the coefficients so that the real and imaginary parts of $f$ separate nicely. In particular, we partition $C \cup S = \{1,\ldots,2n\}$ and then fix  $\eps_{-k} = \eps_k$ for each $k \in C$, and $\eps_{-k} = - \eps_k$ for each $k \in S$. Thus we may write 
\[
 f(x)
 = \eps_0 + 2\sum_{k \in C}\eps_k\cos kx +2i\sum_{k\in S}\eps_k\sin kx. \]
and then define the real and imaginary trigonometric polynomials as 
\begin{equation}\label{eq:poly-form}c(x) =\sum_{k \in C}\eps_k\cos kx \qquad \text{ and } \qquad  s(x) = \sum_{k\in S}\eps_k\sin kx . \end{equation}
While we don't discuss the precise definition of $C,S$ here, it will be important for us that $C \subset [\gamma n]$, where $\gamma$ is a small constant, so that the degree of $c(x)$ is small. 

Thus we construct our function $f$ in two stages. We first construct a cosine polynomial $c(x)$ which is $O(\sqrt{n})$
for all $x$ and satisfies $|c(x)| \geq \delta \sqrt{n}$,  \emph{except} on a set of intervals $\cI = \{ I\}$, which are not too long, well separated and not too numerous. 
In the second, and more challenging step, we shall show that we can construct a sine polynomial $s(x)$ that is $\Omega(n^{1/2})$ on these intervals where the cosine polynomial is small, while still maintaining the upper bound of $O(n^{1/2})$ overall. 

While there are probably many different ways of constructing an appropriate cosine polynomial, we use a deterministic construction based on the Rudin--Shapiro polynomials, mentioned above. Rudin and Shapiro defined their polynomials recursively by defining the pairs of polynomials $P_t,Q_t$ by $P_0(z)=Q_0(z)=1$ and inductively defining 
\[ P_{t+1}(z) =P_t(z) + z^{2^t}Q_t(z),  \qquad \text{ and } \qquad 
  Q_{t+1}(z)= P_t(z) - z^{2^t}Q_t(z).\]
 for each $t \ge 0$.
 
 We construct our cosine polynomial $c(x)$, by using a ``twisted'' version these polynomials. We define  
 \[ c(x) = \Re\big( z^{T} P_t(z) + z^{2T} Q_t(z) \big), \]
where $T \approx \gamma n$ and $t = \log \g n$ so that $\deg(P_t)$ and $\deg(Q_t) \approx \gamma n$. Thus by the boundedness of the Rudin--Shapiro polynomials, we have that 
$|c(x)| = O(\sqrt{n})$. We also have that $|c(x)| \geq \delta \sqrt{n}$ \emph{except} on a collection $\cI$ of intervals that satisfy
\begin{enumerate}
\item $|\cI| = O(\g n) $;
\item Each $I \in \cI$ has $|I| = O(n^{-1})$;
\item For distinct $I,J \in \cI$, we have $d(I,J) = \Omega(n^{-1})$.
\end{enumerate}
Actually, the first condition holds since we arranged for the degree of $c(x)$ to be at most $\g n$. Also note that the second condition holds ``typically'' in the sense that we expect such a polynomial to have derivative $\approx (\g n)^{3/2}$. So we expect it to be within the interval $[-\delta\sqrt{n},\delta \sqrt{n}]$ for time at most $\approx n^{-1}$ (here $\delta \ll \gamma $). We expect the last condition for a similar reason. 

In the second and more challenging part of the proof, we show that if we are given any collection of intervals satisfying the conditions (1)-(3), we can construct a sine polynomial $s(x)$ of the form in \eqref{eq:poly-form} that satisfies $s(x) = O(\sqrt{n})$ everywhere and $|s(x)| \geq \delta \sqrt{n}$ on each interval $I \in \cI$. It's in the construction of $s(x)$ that we use ideas from discrepancy theory. 

Before defining the $s(x)$, we first assign a ``direction'' $\alpha(I) \in \{-1,1\}$ to each bad interval $I\in \cI$, which indicates the sign we want $s(x)$ to have on $I$. (We will describe how we choose these $\alpha_I$ in a moment). We then define, for each $k$, the quantity
\begin{equation}\label{eq:Delta} \Delta(k) =  \sum_{I \in \cI}   \frac{\alpha_I}{|I|}\int_I \sin ks \, ds,\end{equation} which tells us, based on how positive or negative it is, how much we ``prefer'' the choice of $\eps_k = 1$ versus the choice of 
$\eps_k = -1$. Indeed one can think of $\Delta(k)$ as the ``net-progress'' we make towards ``pushing'' our various intervals in the desired directions, when choosing the coefficient $\eps_k$.

To ensure that we carefully ``spread'' each push out over all of the intervals and time steps, we make our first use of discrepancy theory (Theorem~\ref{thm:spencer}) to choose the $\alpha_I$ so that 
\begin{equation}\label{eq:Delta-def} |\Delta(k)| \leq C'|\cI|^{1/2} \leq C(\gamma n)^{1/2} ,\end{equation}
for some absolute constants $C',C>0$.

We now use these quantities $\Delta(k)$ to define a space of \emph{random} sine polynomials. First define the random variables
\[ \hat{\eps}_k \in \{-1,1\} \qquad \text{ by } \qquad \EE\, \hat{\eps}_k = \frac{\Delta(k)}{C(\g n)^{1/2}} ,\] where we are implicitly using \eqref{eq:Delta-def}, to make sure such a random variable exists. We then define the random sine polynomial $\hat{s}(x)$ by 
\[ \hat{s}(x) = \sum_{k \in S} \hat{\eps}_k \sin kx .\] 
Heuristically, the idea is this: with each choice of $\eps_k$, for each $k \in S$, we can increase 
\[ \min_{x \in I} |\hat{s}_{\alpha}(x)|\] by about $\Delta(k)/(C\gamma n)^{1/2} = \Theta((\gamma n)^{-1/2})$, for each $I \in \cI$. Since we have $|S| = \Theta(n)$ values of $k$ to work with, we should (on average at least) push each interval as far as $\geq n^{1/2}/\g^{1/2}$. That is, far enough.

To get some idea why we can indeed guarantee this, we see how orthogonality of the characters, allows us to see that the expected value of $\hat{s}$ is large on \emph{all} of the intervals. Indeed, fix an interval $J \in \cI$ and fix a point $x \in J$. We observe that 
\[ \EE\, \hat{s}(x) = \EE \sum_{k\in S} \eps_k \sin k x =  \frac{1}{C(\gamma n)^{1/2}}\sum_{k\in S} \Delta(k) \sin k x, \]
which gives, expanding the definition of $\Delta(k)$,
\begin{equation}\label{eq:Es-hat} \EE\, \hat{s}(x) =  \frac{1}{C(\gamma n)^{1/2}}\sum_{I \in \cI}   \frac{\alpha_I}{|I|}  \cdot \int_I\,  \sum_{k \in S} (\sin kx )(\sin ks) \, ds .\end{equation}
Now if $d(x,I) < 1/n$ we have the approximate orthogonality relations
\[  |I|^{-1}\int_I \sum_{k \in S} (\sin kx )(\sin ks )  \, ds \approx n \qquad \text{ and } \qquad |I|^{-1}\int_I  \sum_{k \in S} (\sin kx)(\sin ks) \, ds \ll n,  \]
whenever $d(x,I) \gg 1/n$. Thus we have that the sum on the right hand side of \eqref{eq:Es-hat} ``picks out'' the interval $J$. We can therefore conclude that 
\[ \EE\, s(x) \approx \frac{\alpha_Jn}{C(\g n)^{1/2}} = \Theta\big(\sqrt{n/\g}\big) . \]
Thus we have sketched how a sample from $\hat{s}(x)$ behaves correctly on the intervals $I \in \cI$ \emph{on average}. Unfortunately, a \emph{typical} sample from this polynomial will not be enough to push up all the intervals simultaneously. Indeed, the variance is large enough to spoil the value of $\hat{s}(x)$ on many $I \in \cI$. To get beyond this, we appeal to tools in discrepancy theory and, in particular, to a version of the partial colouring lemma mentioned above (Lemma~\ref{lem:recolouring}) due to Lovett and Meka \cite{lovett2015}.

With this we are able to find a (exponentially unlikely) polynomial $s(x)$ with the property that $ |s(x)| \geq \delta n^{1/2}$ for all $x \in \bigcup_{I \in \cI} I$. It is with this polynomial that we can complete the proof.

\subsection{Constructive proofs of Spencer's theorem } 

As we sketched above, Spencer's original proof of Theorem~\ref{thm:spencer} relies fundamentally on an application of the pigeonhole principle and thus, while it does show that a solution $x$ does exist, it gives little guidance on how to find it efficiently. In fact, Spencer conjectured that no efficient algorithm exists.

The first breakthrough was provided by Bansal \cite{bansal2010}, who refuted Spencer's conjecture by providing an efficient algorithm which used a clever random walk which was guided by a semi-definite program, which encoded the current ``state'' of the solution. A few years later a much simpler algorithm was then given by Lovett and Meka \cite{lovett2015}, which has the additional advantage that it did not rely on Spencer's original proof. Because their proof is so simple and elegant, we sketch a proof of their main recoloring step; that is, their analogue of Lemma~\ref{lem:recolouring}.

\begin{lemma}\label{thm:Lovett-Meka} Let $a^{(1)},\ldots, a^{(m)} \in \R^n$ be vectors with $\|a^{(i)}\|_{\infty} \leq 1$ for all $i \in [m]$. Let $x_0 \in [-1,1]^n$ and let $c_1,\ldots,c_m \geq 0$ be such that 
\[ \sum_{j} \exp\big(-c_j^2/16\big) \leq 1/16.\] Then there exists $x \in [-1,1]^n$ so that for all $j \in [m]$, we have 
\[ \big|\la x-x_0, a^{(j)} \ra \big| \leq c_j\sqrt{n} \]
and $|x_i| = 1 $, for at least $n/4$ entries of $x$. 
\end{lemma}

Here we have taken the liberty of removing the quantification that makes it clear that Lemma~\ref{thm:Lovett-Meka} also results in an efficient algorithm, since this is not our focus here. Let us also note that, in contrast to Lemma~\ref{lem:recolouring}, the present lemma fixes the colouring on $n/4$ coordinates and gives a \emph{fractional} weight to all other coordinates. Thus one should think about $x_0$ as the fractional weight ``so far'' in the iteration.

\begin{proof}[Proof sketch]
We imagine two convex bodies. The first is simply the cube 
$[-1,1]^n$, which is defined, of course, by the hyperplanes 
$-1\leq x_i \leq 1$, for $i \in [n]$. 
The second is the convex body defined by the hyperplanes which come from the linear constraints. That is,
\begin{equation}\label{eq:constrains} \big|\la x-x_0, a^{(j)} \ra \big| \leq c_j\sqrt{n}.\end{equation}
We now define the random process $X_t$ as follows. The process $X_t$ starts, at $t=0$, at $x_0$. We then allow $X_t$ to evolve as Brownian motion until it first hits one of the hyperplanes above, at which point it sticks to it and behaves as a Brownian motion \emph{within} the hyperplane. 

Thus our process evolves, wandering as a Brownian motion within a set of hyperplanes, hitting further hyperplanes, and then restricting itself further. Thus $X_t$ is a random process with mean $x_0$ and with a covariance matrix that starts as the identity and which is successively projected onto the hyperplanes that it sticks to.
 
 Using standard martingale concentration estimates, we can say that the large deviations of this process throughout time are no worse than that of the unconditioned Brownian motion. In particular we can see that 
\[ \PP_{X_t}\big( |\la X_t-x_0, a^{(j)} \ra| \leq c_j\sqrt{n} \big) \leq e^{-c_j^2/2}.\]
Thus, by time $t = 1$, the expected number of hyperplanes of type \eqref{eq:constrains} that the process hits is at most $1/16$. Thus it must have hit a good proportion of the hyperplanes of the type $|x_i|=1$, which gives exactly what we want. 
\end{proof}

\subsection{The Koml\'{o}s Conjecture and the Beck-Fiala conjecture}

Before concluding our discussion of discrepancy theory, it is impossible not to mention the beautiful and conjectural extension of Spencer's theorem known as the \emph{Koml\'{o}s conjecture}, which says that one only needs to control the $\ell_2$ norm of the columns of the matrix $A$ to arrive at the same conclusion as Spencer's theorem (the normalization is changed here to match the literature).

\begin{conjecture}
Let $A$ be a $n \times n$ matrix where each column has $\ell_2$-norm at most $1$. Then there exists $x \in \{-1,1\}^n$ so that 
\[ \|Ax\|_{\infty} \leq K, \]
for an absolute constant $K$.
\end{conjecture}

There is also a famous hypergraph colouring ``companion'' to this conjecture, made independently by Beck and Fiala \cite{beck1981}. It is not hard to see that the Koml\'{o}s conjecture implies the Beck-Fiala conjecture.
\begin{conjecture}\label{conj:beck-fiala} Let $\cH$ be a hypergraph on finite ground set where every vertex has degree at most $d$. There exists 
$f : X \rightarrow \{-1,1\}$ so that for all $e \in \cH$ we have 
\begin{equation}\label{eq:beck-fiala} \bigg| \sum_{x \in e} f(x) \bigg| \leq C\sqrt{d},\end{equation}
where $C>0$ can be taken to be a absolute constant. 
\end{conjecture}

Beck and Fiala proved that Conjecture~\ref{conj:beck-fiala} is true if $C\sqrt{d}$ is replaced with $2d-1$. The only unconditional improvement to this bound is by Bukh \cite{bukh2016} who improved it to $2d - \log_{\ast} d$. If we set $|X| = n$, Banaszczyk \cite{banaszczyk1998} proved that one can take $C = O(\sqrt{\log n})$ in \eqref{eq:beck-fiala} by proving one can take $K = O(\sqrt{\log n})$ in the setting of the Koml\'{o}s conjecture. 

These results remained the state of the art for over 25 years, until a very recent and exciting breakthrough of Bansal and Jiang \cite{bansal2025,bansal2025decoupling}. In these papers they prove the Beck--Fiala conjecture in the case $n \leq 2^{c\sqrt{k}}$ and prove a bound of $(\sqrt{k} + \log n)(\log\log n)^{c}$ in general. They also gave the first improvement on Banaszczyk's bound of $K = O(\sqrt{\log n})$ in the setting of the Komlos conjecture, by showing that one may take $K = (\log n)^{1/4+o(1)}$.

\section{Changing the distribution: the semi-random method and sphere packings}\label{sec:processes}

One alternative way of finding a rare object in a probability space is to \emph{change} the underlying distribution that it is sampled from. This can give us a way of naturally ``accessing''  unlikely events. We already saw this kind of idea in action with the constructive proof of Spencer's theorem due to Lovett and Meka, but  perhaps the most classical example comes from the work of Ajtai, Koml\'{o}s, and Szemer\'{e}di \cite{ajtai1980note} on independent sets in triangle free graphs, which was  later refined by Shearer \cite{shearer1983note} to give the following basic and beautiful result. 

For this we recall that an \emph{independent set} in a graph $G$ is a set of pairwise non-adjacent vertices and the \emph{independence number} of a graph $G$, denoted $\alpha(G)$, is the largest independent set in $G$.

\begin{theorem}\label{thm:sheaer} Let $G$ be a triangle-free graph on $n$ vertices with average degree $d$. Then 
\[ \alpha(G) \geq \big(1+o(1)\big) \frac{n \log d}{d} ,\]
where the $o(1)$ term tends to $0$ as $d$ tends to infinity. 
\end{theorem}
This result is easily seen to be sharp, up to a factor of $2+o(1)$ by appropriately modifying a random graph. While this theorem has many applications, perhaps the best known is the following bound on the extreme off-diagonal Ramsey numbers. 

\begin{theorem} We have 
\[ R(3,k) \leq \big(1+o(1)\big) \frac{k^2}{\log k }.\]
\end{theorem}

If one looks to prove Theorem~\ref{thm:sheaer} with a ``direct'' application of the probabilistic method, that is by selecting a set uniformly at random, one is doomed to failure. Indeed in a random graph of average degree $d$ there exponentially few such independent sets among all $k$ sets. To access these sets we instead ``tilt'' our distribution towards independent sets. In fact, there are a couple different ways to do achieve this in practice, but in what follows we outline a heuristic that is behind all of these different approaches.  

\vspace{2mm}

\emph{Heuristic justification of Theorem~\ref{thm:sheaer}.} Suppose that we could define a distribution on independent sets that produced a set that still ``looked random'' apart from the constraint we imposed of being independent. How large could we reasonably expect the size of the independent set we produced? Let's say that our distribution produces an independent set $I$ of $pn$ vertices, for some $p \in (0,1)$. Say a vertex is \emph{open} if $v \not\in I$ and all of its neighbours are not in $I$. The probability that a vertex is left \emph{open} is
\begin{equation}\label{eq:shearer-heuristic} \PP\big( v \text{ open} \big) = \PP\big( (v \cup N(v)) \cap I = \emptyset \big) \approx (1-p)^{d(v)+1},\end{equation} where $d(v) = |N(v)|$ is the degree of $v$ and $N(v)$ is its neighrbouhood. Here we have used our heuristic assumption that $I$ is random-like.
Now intuitively we want to choose $p$ just large enough so that we just start to have 
\begin{equation}\label{eq:p-heuristic} (1-p)^{d} \approx \PP( v \text{ open} ) \ll p .\end{equation} That is, we expect this optimal $p$ to be just at the point where the number of vertices left open is significantly smaller than the number of vertices that we have added so far. Thus, solving for $p$ in the equation \eqref{eq:p-heuristic}, brings us to a heuristic for the maximum density of a ``random-like'' independent set
\[ p = (1+o(1))\frac{\log d}{d}, \]
which exactly matches Sheaer's bound.

Of course this is not a proof at all since we have not provided any distribution that satisfies these conditions. However, there are at least two natural such distributions. The first is to build $I$ by a random greedy process where we remove random vertices one-by-one along with all of their neighbors. This is the idea behind the proofs of Ajtai, Koml\'{o}s and Szemer\'{e}di \cite{ajtai1980note} and the refinement of Shearer. But we also note that a different, more direct proof was given by Shearer~\cite{shearer1995} which uses the hardcore model on $G$ to sample $I$.

\vspace{2mm}

\emph{Towards an optimal version of Shearer's theorem.} We pause to remark that it is a major open problem to determine the correct constant in Theorem~\ref{thm:sheaer}. It is unclear which if either the upper bound or lower bound is sharp. We also mention the beautiful algorithmic problem of finding an independent set in the random regular graph that improves upon Shearer by a constant factor. More precisely, does there exist a randomized polynomial time algorithm that finds an independent set in the random regular graph $G(n,d)$ of size $\geq (1+\eps)n(\log d)/d$, with high probability? In this setting we even know that large independent sets exist and thus should in principle be easier. However it has been shown there are serious obstructions to finding such an algorithm \cite{gamarnik2014limits,rahman2017local}. 

The problem of finding an optimal version of Shearer's theorem is also intimately tied up with the problem of determining the Ramsey numbers $R(3,k)$ and accounts for the missing factor of $2+o(1)$ in this problem \cite{campos2025,hefty2025}.

\subsection{Spherical codes and sphere packing in large dimension}

Recently, in joint work with Marcelo Campos, Mattew Jenssen and Marcus Michelen, we applied this sort of thinking to the classical \emph{sphere packing problem}: What is the maximum proportion of $\R^d$ that can be covered by non-overlapping spheres of volume one?  There is also the closely related question of constructing \emph{spherical codes}: Given an angle $\theta$, what is the maximum proportion of $\mathbb{S}^{d-1}$ that can be covered by non-overlapping spherical caps of radius $\theta$? Let $\theta(d)$ denote this maximum proportion in the sphere packing problem and let $A(\theta,d)$ denote the maximum proportion of $\bS^{d-1}$ in the spherical caps problem. 

Despite the simplicity of these problems, little is known about these fascinating quantities.  The precise value of $\theta(d)$ is only known in dimensions $d\in \{1,2,3,8,24\}$. The case $d=1$ is trivial, the case $d=2$ is classical \cite{thue1911dichteste}, while the remaining known cases are the result of a series of extraordinary breakthroughs: dimension 3 was a landmark achievement of Hales \cite{hales2005proof}, resolving the Kepler conjecture from 1611. Dimensions 8 and 24 were resolved only recently due to the major breakthroughs of Viazovska~\cite{viazovska2017}, in dimension $8$, and then Cohn, Kumar, Miller, Radchenko, and Viazovska~\cite{cohn2016sphere} in dimension $24$. (See~\cite{cohn2016conceptual} for a beautiful exposition of these developments). 

We also recall that the \emph{kissing number} of $\R^d$ corresponds to the special case of the spherical codes problem $A(d,\pi/3)$, although it is more traditionally phrased as the maximum number of unit spheres in $\R^d$ that can be arranged tangent to (or which ``kiss'') a central unit sphere. The only kissing numbers that are known are in dimensions $d \in \{ 1, 2, 3, 4, 8, 24\}$. Similarly only a few cases of optimal spherical codes are known for other $\theta$, for which we refer the reader to~\cite{cohn-website}.

In our work, our focus is on sphere packing and spherical codes in \emph{large dimension}, where the situation remains even more mysterious. A simple argument shows that any \emph{saturated} packing (one in which no additional sphere can be added) has density $\geq 2^{-d}$ and thus 
\[ \theta(d) \geq 2^{-d} . \]
A classical theorem of Minkowski's \cite{Min05} improved upon this bound by a factor of $2+o(1)$. In 1947 Rogers~\cite{rogers1947existence} made the first asymptotically growing improvement to the trivial lower bound showing that 
\[ \theta(d)\geq (\beta+o(1)) d2^{-d}, \] where $\beta=2/e\approx 0.74$. Since the work of Rogers, a number of improvements have been made to the constant factor $\beta$. Davenport and Rogers~\cite{davenport1947hlawka} showed that one can take $\beta=1.68$; Ball~\cite{ball1992lower}, some 45 years later, improved the bound to $\beta=2$; and Vance~\cite{vance2011improved} showed that one can take $\beta=6/e\approx 2.21$ when the dimension $d$ is divisible by $4$. Venkatesh~\cite{venkatesh2012note} showed that one can take $\beta=65963$ and additionally showed that one can obtain an additional $\log \log d$ factor along a sparse sequence of dimensions.  In our paper \cite{CJMS2023sphere-packing}, we go beyond this barrier and improve Minkowski's bound by a factor of $\Omega(d \log d)$ in general dimension. 

\begin{theorem}\label{thm:sphere-packing}
   As $d$ tends to infinity
    \[
    \theta(d) \geq (1-o(1)) \frac{d\log d}{2^{d+1}}\, .
    \] 
\end{theorem}

Recently, this result has been seen a further spectacular improvement by by Klartag \cite{klartag2025}, who used a method reminiscent of Lovett and Meka's proof of Spencer's theorem, to show the following. 
\begin{theorem}\label{thm:klartag} As $d$ tends to infinity
\[ \theta(d) \geq cd^22^{-d},\]
for some $c>0$. \end{theorem}
We discuss this beautiful result further in Section~\ref{sec:Klartag}.

Our method also naturally adapts the setting of spherical codes in large dimension and provides us with an improvement in this setting. To state this result, we let $s_d(\theta)$ denote the normalized spherical volume of a cap of angle $\theta$. In \cite{CJMS2023sphere-packing} we also prove the following.  

\begin{theorem}\label{thm:sphere-codes} If $\theta \in (0,\pi/2)$ and $d$ tends to infinity then
\[ A(d,\theta)\geq (1-o(1)) \frac{d\log d}{2 s_d(\theta)}. \]\end{theorem}
This improved upon the best known bounds due to Fern\'{a}ndez, Kim, Liu and  Pikhurko~\cite{fernandez2025new} who gave a constant factor improvement to bounds of Jenssen, Joos and Perkins~\cite{jenssen2018kissing}.
These bounds were of the type $A(d,\theta) \geq c d /s_d(\theta)$ for some constant $c>0$.

We also note that our results have been adapted further to other settings. Fern\'{a}ndez, Kim, Liu and Pikhurko \cite{fernandez2025-hyperbolic} improved the best bounds for the sphere packing in high dimensional hyperbolic space using this method and Schildkraut \cite{schildkraut2024} has extended this method to show that one can obtain a similar bound for packing balls in an arbitrary norm. 


\vspace{2mm}

\emph{Upper bounds on the sphere packing problem.} Despite this progress, the upper bounds for the sphere packing problem are quite far off the lower bounds, with an exponential gap between the two. The best known upper bounds are of the form 
\[ \theta(d) \le 2^{- (.599\dots +o_d(1)) d},\] which is due to the 1978 work of Kabatjanski\u\i\, and Leven\v ste\u\i n~\cite{kabatiansky1978bounds} and has only been improved by a multiplicative constant factor in the years since by Cohn and Zhao~\cite{cohn2014sphere} and then Sardari and Zargar~\cite{sardari2023new}. It is a beautiful and central problem to improve these bounds further.

 \subsection{Amorphous sphere packings in physics}
One interesting property of the sphere packings behind Theorem~\ref{thm:sphere-packing} is that they are ``random-like''. While essentially all other results focus on \emph{lattice} packings, which are therefore very ``structured'', our packings are essentially as random-like as possible. Such packings are of independent interest in the  physics literature where random sphere packings at a given density are a natural model of physical matter. 

In dimension $3$, for instance, it is believed that random sphere packings transition from appearing ``gas-like'' at low density to ``lattice-like'' at high density, paralleling the phase transition between states of matter. However, rigorously demonstrating that this phase transition occurs remains a beautiful and major open problem in the field (see \cite{lowen2000fun} and the references therein). 

Physicists have also devoted enormous effort to analysing sphere packings in high dimensions, with the aim of providing a more tractable analysis than in low dimensions, and in order to use the powerful machinery of equilibrium statistical physics to generate rich predictions. 
Here, the important qualitative distinction is between sphere packings that are \emph{crystalline}, meaning that  they exhibit long-range ``correlations'', and \emph{amorphous}, meaning they don't have any such correlations.  For example, lattice packings are extreme instances of crystalline packings where the entire structure is determined by a basis.

In their seminal work on applying the replica method to the structure of high-dimensional sphere packings, Parisi and Zamponi~\cite{parisi2010mean,parisi2006amorphous} predicted that the largest density of amorphous packings in $d$ dimensions is 
\[ (1+o(1))(d\log d) 2^{-d},\] that is, a factor of $2$ larger than our lower bound from Theorem~\ref{thm:sphere-packing}. While there is no agreed-upon rigorous definition of ``amorphous,'' it seems likely that any such definition would be satisfied by our construction for Theorem~\ref{thm:sphere-packing}, which enjoys extremely fast decay of correlations.

\subsection{Sketch proof - a graph theoretic reduction}

 To prove Theorem~\ref{thm:sphere-packing} and Theorem~\ref{thm:sphere-codes} we convert the problem into the problem of finding a large independent set in a certain graph. To do this we discretize the space in a natural way. Here we sketch the situation for sphere packings, and note that the case for spherical codes only requires small adjustments. 

To discretize, we simply sample a Poisson point process in a large box $[-L,L]^d$ at intensity 
\[ \lambda  = d^{d/2-o(d)} . \]
We don't worry about the $o(d)$ term, but it is chosen so that, for a typical point in our sample, the next nearest point will be of distance $\gg \log d$. (Some points will have a closer nearest point, but we can simply delete these). Let $X$ be the outcome of this initial discretization step. 

Now a natural graph $G = G_X$ suggests itself. We let $X$ be the vertex set and we define 
\[ x \sim y \qquad \text{whenever} \qquad \|x-y\|_2 < 2r_d , \] where $r_d$ is the radius of a ball of volume one in $\R^d$. That is, $x$ and $y$ are joined by an edge if $B_{r_d}(x) \cap B_{r_d}(y) = \emptyset$. Thus an independent set in $G$ is a sphere packing in the box $[-L,L]^d$.

We now would like to ``lift'' this graph out of its geometric context and think of it only as a graph. But what properties can we hold on to? One obvious one is the degree. We can easily compute the expected degree.
If we fix a point $x \in X$, using the basic properties of Poisson point processes, we can estimate
\[ \EE\, |X \cap B_{2r_d}(x)| = \Vol(B_{2r_d}(x))\lambda = 2^d \lambda =: \Delta .\]
If we were to use this bound along with the trivial bound (mentioned above) 
$\alpha(G) \geq \frac{n+1}{\Delta(G)}$, we can recover the (also trivial) bound $\theta(d) \geq 2^{-d}$. To get beyond this bound we need to use some additional information. Inspired by the theorem of Ajtai, Koml\'{o}s and Szemer\'{e}di (or Theorem~\ref{thm:sheaer}) one might think about focusing on the number of triangles in the graph $G$.
This perspective was taken in \cite{krivelevich2004lower} but only matches the bounds of Rogers and is sharp from this point of view. 

Our new idea is to focus on the maximum \emph{codegree} of our graph, which actually behaves very well in this context. Indeed we can easily compute the co-degree of our graph 
\[ \EE\, |X \cap B_{2r_d}(x) \cap B_{2r_d}(y)| = \Vol(B_{2r_d}(x) \cap B_{2r_d}(y)) \lambda \leq \big( 2^d\lambda \big) e^{-\|x-y\|_2^2/2} \leq \Delta/(\log \Delta)^{\omega(1)}, \]
where in the last inequality we are using that $\|x-y\|_2 \gg \log d$.

The insight here is that we can obtain the same bound as Shearer for graphs that have controlled codegrees. Interestingly, this is also a new result in graph theory. 

\begin{theorem}\label{thm:ind-graph}
Let $G$ be a $n$ vertex graph with $\Delta(G)\leq \Delta$ and $\Delta_2(G)\leq C\Delta (\log \Delta)^{-c}$. Then 
\[\alpha(G)\geq (1-o(1))\frac{n\log \Delta }{\Delta }\,,\]
where $o(1)$ tends to $0$ as $\Delta\to \infty$ and we can take $C = 2^{-7}$ and $c = 7$.
\end{theorem}

\subsection{Sketch proof of Theorem~\ref{thm:ind-graph}} To prove this we use a nibble process as Ajtai, Koml\'os and Szemer\'edi, but our analysis is quite a bit different. We sketch a little to see how the co-degree condition comes naturally into play. As we discussed above, we build our independent set by building it up in pieces. We take our first piece
as $p_1 = \frac{\gamma}{\Delta}$, for some small $\gamma \ll 1$. Let $I_1$ be this $p_1$-random set. Note that since the maximum degree of this graph is $\Delta$, every vertex in $G[I_1]$ will have average degree $\g \ll 1$, and thus $I_1$ is very close to an independent set. Indeed, we can make it independent by throwing away $o(|I_1|)$ vertices.  

We now delete all of $I_1$ and all of the neighbors of $I_1$ from the graph. Define
\[ D_1 = I_1 \cup \bigcup_{x \in I_1} N(x) .\]
which is about $\gamma$ proportion of the vertices of $G$. The key property we would like to maintain is that $D_1$ ``looks like'' a random set of density $\g$ in $G$. If this is possible then we expect that the new maximum degree is about $(1-\gamma)\Delta$ and the new maximum codegree is about $(1-\gamma)\Delta$. Thus we can choose  $p_2 = \g/((1-\gamma)\Delta)$ and then choose $I_2$ to be a $p_2$ random set in the second nibble. Thus we have 
\[ |I_2| \approx p_2 (1-\g)n = \g n/\Delta . \]
More generally, after the $i$th nibble, we will have constructed disjoint sets $I_1, \ldots ,I_i$ with $|I_i| \approx  \g n/\Delta$ and so that $I_1 \cup \cdots \cup I_i$ is independent (after a small amount of clean-up), and the graph remaining after we remove all of the $I_i$ and all vertices adjacent to them has size $(1-\gamma)^in$.  Thus we can continue this process until \[ (1-\g)^in \leq n/\Delta,\]
meaning that we can run the process for $i  \approx (\log \Delta )/\g $ steps. Thus (assuming that we can maintain these properties) we can  construct an independent set of size $\approx (n/\Delta )\log \Delta$.

To make the above story work, the key new idea is in controlling the evolution of the degrees of the vertices. To sketch the idea here, we fix a vertex $x$ and consider $N(x)$ and a stage $i$ of the process. Let us condition on the survival of $x$ into the next process - which means that none of the neighbors of $x$ are selected for $I_i$. Now the size of $N(x)$ is precisely governed by the set 
\[ Y = N(N(x)) \setminus (N(x)\cup \{x\}), \]
the neighbors of the neighbors of $x$, apart from $N(x)\cup \{x\}$ (since $I_i$ will not include vertices of $N(x) \cup \{x\}$). We now run a martingale argument. We iteratively expose each vertex in $Y \cap I_i$. If a vertex $v \in Y$ is included into $I_i$ we then delete all of $N(v) \cap X$ from $X$. Now, to obtain concentration we note that the steps of the martingale are controlled by the sum of the squares of the increments, which due to the double counting inequality
\[ \sum_{y \in Y} |I \cap N(y)|^2 \leq \sum_{y,z \in N(x)} |N(y) \cap N(z)|, \]
are controlled by the co-degrees of the vertices.

\subsection{Klartag's new sphere packing bounds}\label{sec:Klartag}

We now turn to sketch the beautiful new idea of Klartag \cite{klartag2025} that allows one to obtain sphere packings of density $\Omega(d^22^{-d})$. Klartag picks up on an earlier idea of building a packing out of a random lattice. However the novelty in Klartag's proof is that instead of simply selecting the lattice uniformly at random, he cleverly ``guides'' a random process to find a better (and exponentially unlikely) choice.  

The setup is this. We first find a lattice $\Lambda \subset \R^d$ with $\det(\Lambda)=1$ and an ellipsoid $\cE$ of large volume, which is centered at the origin and with $\cE \cap \Lambda = \{0\}$. We then turn this into a sphere packing by applying a linear transformation $T: \R^d \rightarrow \R^d$ with $\det(T)=1$ so that $T(\cE)$ is the Euclidean ball $B$ centered at the origin with $\Vol(B) = \Vol(\cE)$. Note that $T(\Lambda)$ is a new lattice with determinant one. Thus if we place a copy of the dilated ball $B/2$ at each lattice point of $T(\Lambda)$ we obtain a sphere packing of identical balls with density $\Vol(\cE)2^{-d}$.

\begin{proof}[Proof sketch of Theorem~\ref{thm:klartag}]
By the discussion above, we see that the problem reduces to the problem of finding a lattice $\Lambda$ with $\det(\Lambda)=1$ and a centrally symmetric ellipsoid of volume $\Omega(d^2)$ that contains no lattice points of $\Lambda$, apart from the origin. In a first step, we choose $\Lambda$ to be a random lattice with determinant $1$. While we won't say anything technically about how to work with these lattices here, it is enough to say that this lattice looks ``locally'' like a Poisson point process with intensity $1$.

We then grow an ellipsoid $\cE_t$ in a manner analogous to the proof of Lemma~\ref{thm:Lovett-Meka}, although here we are working in the space of ellipsoids. Let $\cE_0$ be a euclidean ball which is small enough to ensure that $\Lambda \cap \cE_0 = \{0\}$. We then randomly ``evolve'' this ellipsoid as time proceeds. As soon as this ellipsoid hits a lattice point, it sticks to it and evolves further, keeping this point on its boundary. Indeed, we may describe the ellipsoid $\cE_t$ as
\[ \cE_t = \big\{ x \in \R^d : \la x , A_tx \ra \leq 1 \},\]
where $A_t$ is a positive definite matrix. Thus hitting a point $y \in \Lambda$, precisely introduces the linear constraint $\la y , A_t y \ra = 1$ on $A_t$. Since the dimension of the space of such positive semi-definite ellipsoids is $\approx d^2/2$, we expect that the process runs until the ellipsoid has $\approx d^2/2$ points of $\Lambda$ on its boundary.  

Thus we can heuristically argue about the volume of the final ellipse. If $|\cE_{T} \cap \Lambda| \approx d^2$, one can use the fact that the random lattice $\Lambda$ locally looks like a Poisson point process to see that 
\[ d^2 \approx \EE\, |\cE_{T} \cap \Lambda|   \approx \Vol(\cE_T),\]
as desired. \end{proof}

\section{Random matrix theory at exponentially small scales}
 We now turn to discuss phenomena in random matrix theory that occur at exponentially small scales. Here we focus on the singularity probability of a random symmetric matrix. 

 Let $B_n$ be a random $n \times n$ matrix whose entries are chosen independently and uniformly from $\{-1,1\}$. It is an old problem, likely stemming from multiple origins, to determine the probability that $B_n$ is singular. While a moment's thought reveals the lower bound of $(1+o(1))2n^2 2^{-n}$, the probability that two rows or columns are equal up to sign,
    establishing the corresponding \emph{upper bound} remains an extremely challenging open problem. Indeed, it is widely believed that 
    \begin{equation} \label{eq:asym-conj} \P(\det(B_n) = 0) = (1+o(1))2n^2 2^{-n}. \end{equation}
    While this precise asymptotic has so far eluded researchers, some stunning advances have been made on this fascinating problem. The first steps were taken by the pioneering work of
    Koml\'os \cite{komlos} in the 1960s, who showed that the singularity probability is $O(n^{-1/2})$.
    
    Nearly thirty years later, Kahn, Koml\'os and Szemer\'edi \cite{KKS}, in a remarkable paper, showed that the singularity probability is exponentially small. 
    At the heart of their paper is an ingenious argument with the Fourier transform that allows them to give vastly more efficient descriptions of ``structured'' subspaces of $\R^n$ that are spanned by $\{-1,1\}$-vectors. Their method was then developed by Tao and Vu \cite{TV-RSA,TV-JAMS} who showed a bound of $(3/4)^{n+o(n)}$, by providing a link between the ideas of \cite{KKS} and the structure of set addition and, in particular, Freiman's theorem (see \cite{tao-vu-book}). This trajectory was then developed further by Bourgain, Vu and Wood \cite{BVW}, who proved a bound of $2^{-n/2 + o(n)}$, and by Tao and Vu \cite{tao-vu-ILO}, who pioneered the development of ``inverse Littlewood-Offord theory'', which we disucss below. 
    
    In 2007, Rudelson and Vershynin, in an important and influential paper \cite{RV}, gave a different proof of the exponential upper bound on the singularity probability of $B_n$. The key idea was to construct efficient $\eps$-nets for points on the sphere that have special anti-concentration properties and are thus more likely to be in the kernel of $B_n$. This then led them to prove an elegant inverse Littlewood-Offord type result, inspired by \cite{tao-vu-ILO}, in a geometric setting.
    
    This perspective was then developed further in the breakthrough work of Tikhomirov \cite{Tikhomirov}, who proved 
    \begin{equation}\label{eq:tik} \P(\det (B_n) = 0) = 2^{-n + o(n)},\end{equation} thereby proving the conjectured upper bound, up to subexponential factors. One of the key innovations in \cite{Tikhomirov} was to
    adopt a probabilistic viewpoint on such Littlewood--Offord questions,  a topic which we elaborate on in Section~\ref{sec:LwO}
    
    
    We remark that another pair of advances was made by Jain, Sah and Sawhney~\cite{jain2021}, following the work of Litvak and Tikhomirov \cite{Lit-Tik}, who proved the natural analogue of \eqref{eq:asym-conj} for random matrices with lopsided distributions. In the case of $\{-1,1\}$-matrices, however, the problem of establishing \eqref{eq:asym-conj} perhaps remains as the central open problem in the area. 

\vspace{2mm}

\emph{Singularity of random symmetric matrices.} We now turn to discuss the singularity problem for random \emph{symmetric} matrices, which has proven to be more challenging still. The study of random symmetric matrices goes back to the pioneering work of Wigner in the 1950s (sometimes such random matrices are called \emph{Wigner} matrices) who studied the typical eigenvalue distribution of these matrices, showing that they follow the so called ``semi-circular law''.

 Let $A_n$ be drawn uniformly at random among all $n \times n$ symmetric matrices with entries in $\{-1,1\}$. Again we have the lower bound 
\begin{equation}\label{eq:sing-lb} \PP\big( \det(A_n) = 0 \big)  \geq 2^{- n + o(n)},\end{equation}
by considering the probability that two rows are equal up to sign.  Costello, Tao and Vu \cite{costello-tao-vu} were the first to show that $A_n$ is non-singular with high probability. That is, 
\[ \PP\big( \det(A_n) = 0 \big) = o(1),\]
with a precise error term of $O(n^{-1/4})$. Since, this result has seen a sequence of improvements: A bound of $N^{-\omega(1)}$ was proved by Nguyen \cite{nguyen-singularity}, 
a bound of the form $\exp(-n^{c})$ was proved by Vershynin \cite{vershynin-invertibility}, which was in turn improved by Ferber and Jain \cite{ferber-jain} based on the results of Ferber, Jain, Luh and Samotij \cite{ferber2019counting}. In a similar spirit, Campos, Mattos, Morris and Morrison \cite{CMMM} then improved this bound to $\exp(-c n^{1/2})$  by proving a``rough'' inverse Littlewood-Offord theorem, inspired by the theory of hypergraph containers. This bound was then improved by Jain, Sah and Sawhney \cite{JSS-symmetric}, who improved the exponent to $-c n^{1/2} \log^{1/4} n$, and, simultaneously, by the author to $-c (n \log n)^{1/2}$ in joint work with Campos, Jenssen and Michelen \cite{RSM1}. 
    
As might be suggested by the ``convergence'' of these results to the exponent $-c (n \log n)^{1/2}$, a natural barrier lurks exactly at this point. In fact, in \cite{CMMM} authors showed that if one wants to get beyond bounds at this probability scale, one needs use ``reuse'' randomness in the top half of the matrix (which is of course independent) in the most difficult part of the proof. Rather one needs to directly deal with the complicated dependencies that are present in a random symmetric matrix.

In recent work the author, in joint work with Campos, Jenssen and Michelen, managed to get around this obstacle and prove an exponential upper bound, thus matching \eqref{eq:sing-lb} up to base of the exponent. 
\begin{theorem}\label{thm:singularity} Let $A_n$ be drawn uniformly at random from the $n\times n$ symmetric matrices with entries in $\{-1,1\}$. Then 
\[ \PP\big( \det(A_n) = 0 \big) \leq e^{-cn} ,\]
where $c>0$ is an absolute constant.
\end{theorem}

In what follows we discuss some of the techniques that are behind this result. This will allow us to touch on some of the exciting ideas that have been developed in this area. 

\vspace{2mm}

\emph{Least singular value, clustering and repulsion of eigenvalues.} The singularity problem is related to several other phenomena regarding the spectrum of the matrix $A_n$, the most natural being the extreme behavior of the least singular value. Recall that if $M$ is an $n\times n$ matrix the least singular value is $\sigma_{\min}(M) = \min_{x \in \bS^{n-1}} \|Mx\|_2$.
The study of this quantity in random matrices was first initiated by Goldstine and von Neumann~\cite{goldstine1951} in the 1950s and has undergone intense study in the intervening years, partly in its own right, but also because of its link with spectral laws of random matrices \cite{TV10d,RT19,SSS23b,SSS23} and the smoothed analysis of algorithms \cite{Spielman2002}.

A key guiding conjecture is due to Spielmen and Teng \cite{Spielman2002}, who suggested that in the case of iid Bernoulli random matrices $B_n$ we have 
\begin{equation}\label{eq:spielman-teng}
 \PP( \s_{\min}(B_n)  \leq \eps n^{-1/2} ) \leq \eps + 2e^{-cn},
\end{equation}
for all $\eps>0$.

A key breakthrough on this problem was made by Rudelson \cite{rudelson-annals} which inspired a sequence of further papers \cite{tao-vu-ILO, RV09}, 
cumulating in the influential papers of Rudelson and Vershynin \cite{RV}, who proved \eqref{eq:spielman-teng} up to a constant factors, and Tao and Vu~\cite{TV10dist} who proved \eqref{eq:spielman-teng} in the case of $\eps> n^{-c}$. Recently, in joint work with Sah and Sawheny \cite{sah2025}, the author proved \eqref{eq:spielman-teng}, up to a $1+o(1)$ factor. 

This question has also been intensely studied in the case of random symmetric matrices. In this case we have the additional interpretation that $\sigma_{\min}(A) = \min_i |\lambda_i|$, where the minimum is over all the eigenvalues of $A$. After many partial advances \cite{nguyen-singularity,nguyen-least-singular-value, vershynin-invertibility, JSS-symmetric} in the paper \cite{RSM3}, we determined the optimal probabilities of having a small least singular value for such random symmetric matrices. We showed that for all $\eps >0$
\begin{equation}\label{eq:CJMS-lsv} \PP\big( \sigma_n(A_n) \leq \eps n^{-1/2} \big) \leq C\eps + e^{-cn},\end{equation} where $C,c>0$ are absolute constants. 

One can also apply these results to understand the clustering of the spectrum more generally. Indeed we can apply a version of \eqref{eq:CJMS-lsv} to the matrix $A_n-\lambda I$, for any $-(2-\delta)\sqrt{n} \leq \lambda \leq (2-\delta)\sqrt{n}$, to bound the probability of the event $\min_i|\lambda_i -\lambda| \leq \eps n^{1/2}$.

Allowing ourselves to speak somewhat informally, the form of this result, with two different terms in the bound, reflects two very different phenomena that the least singular value can have. If $\eps \gg e^{-cn}$, for some $c$, the most likely way to have $\sigma_{\min}(A_n) < \eps n^{-1/2}$ comes from the event of a single ``random-like'' direction being hit very weakly by $A_n$. On the other hand if $\eps$ is a small enough exponential, the most likely way to have $\sigma_{\min}(A_n) < \eps n^{-1/2}$ comes from the matrix just being singular, which (conjecturally) should come from the very structured every of having two rows or columns equal. 

These sort of problems are also related to further questions of clustering and repulsion of eigenvalues and we refer the reader to \cite{nguyen2017, nguyen2018}.

\vspace{2mm}

\emph{Anti-concentration of the determinant and permanent.} Before we turn to discuss techniques, we turn to highlight two further problems in this area. The first concerns the anti-concentration of the determinant. Concretely, what is $\PP( \det(B_n) = 1 )$? (or similarly for the symmetric model $A_n$). Here the above proofs for singularity, immediately give that this probability is exponentially small, but conjecturally should be much smaller, most likely of the form $n^{-cn}$. Indeed it seems to be an major step to prove that this probability is $e^{-\omega(n)}$.

We also highlight the related problem on the \emph{permanent} of a random matrix. In this case there is no natural geometric interpretation, forcing us reason by other means. Tao and Vu \cite{tao2009} showed that $\Per(B_n) = 0 $, with probability $O(n^{-c})$ which was the state of the art until a very recent breakthrough on this problem by Hunter, Kwan and Sauermann \cite{hunter2025} who have showed an exponential upper bound. Similar to the previous question, it seems that a bound of the type $n^{-cn}$ should be the truth.

\subsection{Littlewood--Offord theory }\label{sec:LwO} Littlewood--Offord theory is a classical subject that has undergone intense development in recent years as it has become interwoven with the methods used in random matrix theory. The main object of study here is the \emph{concentration function} $\rho_{\eps}(v)$, where $v \in \R^n$ and $\eps >0$, which is defined by  
\[ \rho(v) = \max_{ b } \, \PP\big(  \big| \la X, v \ra - b \big| < \eps \big), \] where $X \in \{-1,1\}^n$ is sampled uniformly at random and the maximum is over all $b \in \R$. (Actually much more general distribution for the $X$ are considered, but we limit ourselves to uniform $X$ here). 

To get a feel about how this is immediately connected to the problems above, we consider the singularity problem for \emph{iid} matrices (i.e. not symmetric). As above, we let $B_n$ be an $n\times n$ matrix with all entries independent and uniform in $\{-1,1\}$. We now expose
the first $n-1$ rows of the matrix $B_n$ and define $\cE$ to be the event that the first $n-1$ rows have full rank. Thus, if we let $X \in \{-1,1\}^n$ be the last row of $B_n$ (and so far unexposed), we have that $B_n$ is singular if and only if $\la X,v \ra = 0$, on the event $\cE$. Thus 
\begin{equation}\label{eq:LwO-and-sing} \PP( \det(B_n) = 0 ) \leq  \PP( \la X, v \ra = 0 \text{ and }  \cE ) + \PP(\cE^c) \leq \EE_{v}\, \rho_0(v) + \PP(\cE^c), \end{equation}
where the expectation is over vectors $v$ that occur as normals to the subspace defined by the first $n-1$ rows. While the probability $\PP(\cE^c)$ can be taken care of by induction, or otherwise, the main difficulty is in dealing with $\rho_0(v)$.

While one expects a somewhat typical vector $v$ to have $\rho_0(v) \approx 2^{-n}$ (further supporting our intuition for \eqref{eq:sing-lb}), there exist $v$ for which $\rho_0(v)$ can be as large as $1/2$, for example $v = (1,-1,0,\ldots,0)$. Moreover, anything in between these two extremes is possible. Thus, the central challenge in estimating the singularity probability is to show that it is unlikely that a vector $v$ with large $\rho_0(v)$ will be orthogonal to the first $n-1$ rows. Thus we are led to understand the concentration function $\rho_0(v)$, as $v$ varies over all possible normals.

\vspace{2mm}

\emph{Classical theory.} Interestingly, the study of $\rho_0(v)$ long pre-dates the study of random matrices, going back to the work of Littlewood and Offord~\cite{LwO-1,LwO-3} in the 1930s on the zeros of random polynomials. (Actually we already mentioned these papers in our discussion of flat Littlewood polynomials.) In 1945 Erd\H{o}s \cite{erdos-LwO} proved what is perhaps the subject's first flagship result, showing that if $v \in \R^n$ has all non-zero coordinates then 
\[ \rho_0(v) \leq \rho((1,\ldots, 1)) = O(n^{-1/2}).\] This was then developed by Szemer\'{e}di and S\'{a}rk\"{o}zy~\cite{SzemerediUNDSarkozy}, Stanley \cite{stanley} and many others \cite{sali,griggs,kleitman}. These early results provide us with a beautiful combinatorial perspective on the problem, but most important for us is the pioneering work of Hal\'{a}sz \cite{halasz} who made an important connection with the Fourier transform, thus giving us a different \emph{analytic} perspective on the problem. 

\vspace{2mm}

\emph{Inverse Littlewood Offord theory.} More recently the question has been turned on its head by Tao and Vu \cite{tao-vu-ILO}, who pioneered the study of ``inverse'' Littlewood-Offord theory. They suggested the following  ``meta-conjecture'' that has guided much subsequent research.

\begin{center}
If $\rho_0(v)$ is ``large'' then $v$ must exhibit arithmetic structure.
\end{center}
\noindent This ``meta-conjecture'' has been addressed in the work of Tao and Vu~\cite{tao-vu-ILO, tao2010sharp}, and Nguyen and Vu \cite{nguyen-vu-1, nguyen-vu-2} who proved that if $v$ is such that $\rho(v) > n^{-C}$ then $O(n^{1-\eps})$ of the coordinates $v_i$ of $v$ can be efficiently covered with a generalized arithmetic progression of rank $r = O_{\eps,C}(1)$. 

While these results provide a very satisfying picture in the range $\rho_0(v) > n^{-C}$, they begin to break down when $\rho(v) = n^{-\omega(1)}$ and are therefore of limited use at exponential probability scales. More recently these ideas have been extended to give structural information about $v$ when $\rho_0(v)$ as small as $\exp(-c\sqrt{n\log n})$, but these results are of a somewhat different nature \cite{ferber2019counting,CMMM,RSM1}.

Perhaps the most relevant among these results for our discussion concerning Theorem~\ref{thm:singularity} is the inverse Littlewood-Offord theorem of Rudelson and Vershynin~\cite{RV} that allows one to control $\rho(v)$ in terms of a related quantity known as the ``least common denominator`` (LCD) of the vector $v$. This result gives relatively weak information about $v$, relative to the results mentioned above, however is effective and very useful all the way down to exponential scales. As this quantity will pop up in our own work, we postpone the discussion of this to Section~\ref{sec:approx-neg-dependence}. 

\vspace{2mm}


\emph{Typical Littlewood-Offord theory.} More recently, with the breakthrough work of Tikhomirov \cite{Tikhomirov} (discussed at \eqref{eq:tik}), a fresh perspective has been shed on \eqref{eq:LwO-and-sing}. Instead of trying to stratify the different behaviors of $\rho_0(v)$ with different ``inverse'' theorems, he directly studies the distribution of $\rho_0(v)$ as a random variable. More precisely, he considers $\rho_{\eps}(v)$ at each fixed scale $\eps  > 2^{-n +o(n)}$, where $v \sim \mathbb{S}^{n-1}$ is chosen uniformly at random from the unit sphere. Now for such random $v$ and $\eps > 2^{-n + o(n)}$ one has
\[ \EE_v\, \rho_{\eps}(v) = \Theta( \eps ). \]
The technical heart of the work \cite{Tikhomirov} is the following tail-estimate on the distribution of $\rho_0(v)$;  
\[ \PP_v\big( \rho_{\eps}(v) \geq L\eps  \big) \leq L^{-\omega(n) } \]
for appropriately large (but fixed) $L$. In our work on the singularity probability, we also take a probabilistic perspective but employ a completely different set of techniques to understand these sorts of tail events, a topic we discuss in Section~\ref{sec:rank-splitting}.

\subsection{Approximate negative correlation} \label{sec:approx-neg-dependence} One of the new ingredients introduced for Theorem~\ref{thm:singularity} is an
``approximate negative correlation'' inequality for linear events. We first discuss this result in its own right and then sketch how it fits into place in the proof of Theorem~\ref{thm:singularity} in in Section~\ref{sec:rank-splitting}. 

We say that two events $A,B$ are \emph{negatively correlated} if 
\[ \PP(A\cap B) \leq \PP(A) \PP(B).\]
In what follows we let $\eps >0$ and let $X \in \{-1,1\}^n$ be a uniform random vector. Here we are interested in ``linear'' events of the shape \begin{equation} \label{eq:small-eps-ball}|\la X, v\ra| \leq \eps . \end{equation}
The result we discuss here shows the approximate negative dependence of the event \eqref{eq:small-eps-ball}, for all $\eps > e^{-cn}$, from the intersection of events 
\begin{equation} \label{eq:intro-two-events} |\la X, w_1 \ra| \leq \beta,\, |\la X, w_2 \ra| \leq \beta, \ldots  \quad , |\la X, w_k \ra| \leq \beta,
\end{equation}
where $\beta>0$ is small but fixed and $w_1,\ldots,w_k$ are orthonormal  vectors with $k\leq cn$. Crucially, in this statement we don't assume anything about the structure of the vectors $w_i$
and allow the dimension of the space they span to be as large as $\Theta(n)$. Our main negative dependence result says something of the general shape 
\begin{equation}\label{eq:small-ball-splitting}
\PP_X\bigg(\big\{ |\la X, v\ra| \leq \eps \big\} \cap \bigcap_{i=1}^k \big\{ |\la X,  w_i \ra| \leq \beta \big\}\bigg) 
\leq \PP_X\big( |\la X, v \ra| \leq \eps  \big)\PP_X\bigg( \bigcap_{i=1}^k 
\big\{ |\la X, w_i \ra| \leq \beta \big\}\bigg),\end{equation}
although in an approximate form. 

To state this result properly, we use the notion of the ``Least Common Denominator'' of a vector $v$, introduced by Rudelson and Vershynin~\cite{RV} and mentioned above in our discussion of Littlewood--Offord theory. For $\alpha \in (0,1)$, we define
the \emph{least common denominator} of a vector $v \in \R^n$ to be\footnote{For a $A \subseteq \R^n$ and $x \in \R^n$ then define $d(x,A) := \inf_{a\in A}\{ \|x-a\|_2 \}$, as is standard.}
\[D_\alpha(v)= \inf\left\{\phi>0:~d(\phi \cdot v, \Z^n \setminus \{0\}) \leq \sqrt{\alpha n} \right\}.\] 
 That is, this quantity is the smallest multiple of the vector $v$ for which it is ``close'' to the integer lattice $\mathbb{Z}^n$ (the calue of. Rudelson and Vershynin showed that $(D_{\alpha}(v))^{-1}$ behaves quite a bit like $\PP( |\la X,V \ra| \leq \eps ) $.

 \begin{theorem}\label{lem:RVsimple}
    For $d\in \N$, $\alpha\in (0,1)$ and $\eps >0$, let $v\in \bS^{n-1}$ satisfy $D_{\alpha}(v)> c/\eps$. If $X \sim \{-1,1\}^d$ is uniform then 
    \[
    \P\big( |\langle X , v \rangle|\leq \eps  \big)\leq C \eps + e^{-c\alpha n}\, ,
    \] where $C, c >0$ are absolute constants. \end{theorem}

Our first main negative dependence proves an approximate version of \eqref{eq:small-ball-splitting} with $(D(v))^{-1}$ (a slightly better behaved quantity) as a proxy for $\PP( \la X,v \ra \leq \eps)$. The following is a formal statement. 

\begin{theorem}\label{thm:decouple 1} Let $n \in \N$, $\alpha \in (0,1)$, $0\leq k \leq  \alpha \beta n$ and $\eps \geq \exp(-\alpha \beta n)$. Let $v \in \bS^{n-1}$ and let $w_1,\ldots,w_k \in \bS^{n-1}$ be orthogonal. If $X \in \{-1,1 \}^n$ is a uniform random vector and $D_{\alpha}(v) > 16/\eps$ then 
\begin{equation}\label{eq:thminvLwO} \PP_X\bigg( |\la X, v \ra| \leq \eps \, \text{ and }\,  \bigcap_{i=1}^k 
\big\{ |\la X, w_i \ra| \leq \beta \big\} \bigg) < C \eps e^{- c k}, \end{equation}
where $C,c>0$ are absolute constants. 
\end{theorem}

Our proof in \cite{RSM2} of Theorem~\ref{thm:decouple 1} is a delicate argument with the $O(n)$-dimensional Fourier transform. As the proof of this somewhat different from the other results in this survey we don't elaborate on it further. 

\subsection{Sketch of the proof of Theorem~\ref{thm:singularity}} \label{sec:rank-splitting} Now that we have motivated a few of the tools behind Theorem~\ref{thm:singularity}, we now turn to sketch its proof. In analogy with our discussion in Section~\ref{sec:LwO}, we study the ``$n$-dimensional concentration function'' which we define to be 
\begin{equation}\label{eq:feps-def} f_{\eps}(v) = \PP_{A_n}\big( \|A_nv\|_2 \leq \eps n^{1/2} \big),\end{equation} where $v \in \bS^{n-1}$, $\eps>0$ and $A_n$ is a uniformly drawn from $n\times n$ symmetric matrices with entries in $\{-1,1\}$.

Intuitively speaking, we expect $f_{\eps}(v)$ to be ``large'' for directions $v$ that are more likely to appear in the kernel of $A_n$ and therefore, to understand the singularity probability
of a matrix $A_n$, it is essential to understand the \emph{upper tails} of the random variable $f_{\eps}(v)$ when $v \sim \bS^{n-1}$ is sampled uniformly at random. As we discussed above, this probabilistic interpretation of singularity was pioneered by Tikhomirov \cite{Tikhomirov} and is a convenient perspective for us to adopt here, although our techniques are quite different. 

Moreover, if we want to prove \emph{exponential} bounds on the singularity 
probability, we need to control this function for $\eps $ as small as $e^{-cn}$. For technical reasons we also need to restrict ourselves to vectors on the sphere that have controlled infinity norm. We call this (vast majority) of the $n$-dimensional sphere $\bS_{\flat}^{n-1}$. Central to our proof is a large deviation estimate of the following type.

\begin{theorem}\label{thm:concentrationOff} For $L > 1 $ and $e^{-cn} \leq \eps \ll 1$ we have that 
\[ \PP_{v}\big( f_{\eps}(v) \geq (L\eps)^n \big) \leq (cL)^{-2n}, \]
where $v \sim \bS_{\flat}^{n-1}$ is sampled uniformly at random and $c>0$ is an absolute constant.
\end{theorem}

In what follows, we sketch the proof of a weaker form of Theorem~\ref{thm:concentrationOff} where we prove a bound of the type $(cL)^{-n}$ in place of $(cL)^{-2n}$. This former bound is too weak for our purposes, but most of the main ideas are contained in its proof. Indeed to prove this, it is enough to show
\begin{equation}\label{eq:ex-f-eps} \EE_{v}\, f_{\eps}(v)= \EE_{v}\, \PP_{A_n}\big( \|A_n v\|_2 \leq \eps n^{1/2} \big)  \leq (C\eps)^{n} \end{equation}
and can then apply Markov's inequality to finish. 

To this end, our first step is to break up the sphere based on the set of coordinates are well behaved. Indeed by now standard methods, we can assume that $v_i = \Theta(n^{-1/2})$ for $d = cn$ values of $i$. By union bounding over all such choices for these coordinates, it is enough to assume $v_i = \Theta(n^{-1/2})$ for all $i \in [d]$.

We then show that we can ``replace'' the matrix $A_n$ (in the definition of $f_{\eps}$ in \eqref{eq:feps-def}) with a random matrix $M_n$ that has many of the entries zeroed out. This will allow us to focus on the well behaved part of $v_i$ and additionally untangle some of the more subtle and complicated dependencies. Indeed we show, by an appropriate Fourier argument, that 
\begin{equation}\label{eq:replacement} f_{\eps}(v) \leq C^n \cdot \PP_{A_n} \big( \|M_nv \|_2 \leq \eps  n^{1/2} \big) ,\end{equation} for some $C>1$, where $M_n$ is the random matrix defined by\footnote{For a vector $v \in \R^n$ and $S \subseteq [n]$, we use the notation $v_S := (v)_{i\in S}$
 and for a $n \times m$ matrix $A$, and $R \subseteq [m]$, we use the notation $A_{S \times R}$ for the $|S| \times |R|$ matrix $(A_{i,j})_{i \in S, j \in R}$. }
\begin{align}\label{eq:Mdef}
 M =  
\begin{bmatrix}
{\bf 0 }_{[d] \times [d]} & H^T \\
H & { \bf 0}_{[d+1,n] \times [d+1,n]}  
\end{bmatrix}. \end{align}
Here $H$ is a $(n-d) \times d$ random matrix with iid\ entries that are $\mu$-lazy, meaning that $(H)_{i,j} = 0$ with probability $1-\mu$ and $(H)_{i,j} = \pm 1 $ with probability $\mu/2$ for some appropriately small $\mu$.

We now use this special form of $M$ to break up the event 
$ \| Mv \|_2 \leq \eps n^{1/2}.$ Indeed, we write 
\[
Mv= \begin{bmatrix}
H v_{[d]}  \\
 H^T v_{[d+1, n]}

\end{bmatrix}
\] 
and so we need only to control the intersection of the events  
\[ \big\|H v_{[d]}\big\|_2  \leq \eps n^{1/2} \qquad \text{ and } \qquad \big\| H^T v_{[d+1, n]}\big\|_2 \leq \eps n^{1/2}.\] Note that if we simply \emph{ignore} the second event and work only with the first, we land in a situation very similar to previous works; where half of the matrix is neglected entirely and we are thus limited by the $(n\log n)^{1/2}$ obstruction, discussed above Theorem~\ref{thm:singularity}. To overcome this barrier, we crucially need to control these two events \emph{simultaneously}. The key idea is to use the randomness in $H$ to control the event $ \| H v_{[d]} \|_2 \leq \eps n^{1/2}$ and we use the randomness in the selection of 
$v \in \bS_{\flat}^{n-1}$ to control the event $ \| H^T v_{[d+1, n]} \|_2 \leq \eps n^{1/2}$.

For this, we partition the outcomes in $H$, 
based on a robust notion of rank (hence the name ``rank-splitting''). We define the event $\cE_k$ to be the event that $H$ has at $k$ ``unhealthy'' singular values,
\[ \cE_k = \big\lbrace H : \s_{d-k}(H)\geq c\sqrt{n} \, \text{ and } \, \s_{d-k+1}(H)<  c\sqrt{n} \big\rbrace\,,\] 
where $\s_1(H)\geq \cdots \geq \s_d(H)$ denote the singular values of $H$. We then bound 
\[ \PP_{M}\big( \| Mv\|_2 \leq \eps n^{1/2} \big) \] above by (only using the randomness in $M$, for the moment)
\begin{equation}\label{eq:ranksplit}
\sum_{k=0}^d  \PP_{H}\Big( \|H^T v_{[d+1, n]}\|_2 \leq \eps n^{1/2}\, \big|\, \|H v_{[d]}\|_2 \leq \eps n^{1/2}\, \wedge\, \cE_k \Big) \cdot
\PP_{H}\Big(\|H v_{[d]}\|_2 \leq \eps n^{1/2} \, \wedge \, \cE_k \Big)\, .\end{equation}
We now see the link with our ``approximate negative dependence theorem'', which we discussed in section~\ref{sec:approx-neg-dependence} which we use (after a good deal of preparation) to bound the quantity
\[ \PP_{H}\big( \|H v_{[d]}\|_2 \leq \eps \sqrt{n} \, \wedge \,\cE_k \big).\] Indeed, after ``tensorizing'' Theorem~\ref{thm:decouple 1} and approximating these objects with appropriate nets, we are able to conclude that 
\[ \PP_{H}\big( \| H v_{[d]} \|_2 \leq \eps \sqrt{n},\, \wedge\, \cE_k \big)\leq (C\eps e^{-c k})^{n-d},\]
unless $v_{[d]}$ is ``structured'', in which case we do something different (and substantially easier). Thus, for all non-structured $v$, we have \eqref{eq:ranksplit} is at most something of the form 
\begin{equation}\label{eq:ranksplit2} 
(C\eps)^{n-d} \sum_{k=0}^d e^{-c k(n-d)} \cdot \PP_{H}\Big( \|H^T v_{[d+1, n]}\|_2 \leq \eps n^{1/2}\, \big\vert\, \|H v_{[d]}\|_2 \leq \eps n^{1/2} \, \wedge \, \cE_k\Big) .\end{equation}
Up to this point, we have not appealed to the randomness in the choice of $v \in \bS_{\flat}^{n-1}$, beyond imposing that $v$ is non-structured. We now introduce the randomness in $v$ by taking an expectation in $v$, bringing us back to our goal of bounding \eqref{eq:ex-f-eps}. Taking expectations and swapping the order of the expectations above, tells us that \eqref{eq:ex-f-eps} is at most 
\begin{equation}\label{eq:ranksplit3}
(C\eps)^{n-d} \sum_{k=0}^d e^{-c k(n-d)} \cdot \EE_{H}\, \PP_{v}\Big( \|H^T v_{[d+1, n]}\|_2 \leq \eps n^{1/2} \Big) \one\big( H \in \cE_k \big) .\end{equation}
We then deal with this inner probability by considering a \emph{fixed} $H \in \cE$. Here one can show 
\begin{equation}\label{eq:expvsketch}
  \PP_{v_{[d+1,n]}}\Big( \|H^T v_{[d+1, n]}\|_2 \leq \eps n^{1/2} \Big) \leq (C\eps)^{d-k}\, .\end{equation} 
  Indeed, intuitively this is clear, since $H^{T}$ has at most $k$ small singular directions, so $v_{[d+1,n]}$ must be ``nearly orthogonal'' to the $d-k$ singular directions of $H$. At this point one might be slightly worried that we only have another hard high-dimensional Littlewood-Offord problem on our hands. However, the big advantage here is that $v_{[d+1,n]}$ is a \emph{continuous} random variable and thus its analysis is vastly easier.

Now stringing together \eqref{eq:replacement},\eqref{eq:expvsketch} and \eqref{eq:ranksplit3} (and using that $\eps > e^{-cn}$) we arrive at our goal of showing that
\[
\EE_v\, f_{\eps}(v) \leq C^n\EE_v\, \PP_{M}\big(\| Mv\|_2 \leq \eps n^{1/2} \big) \leq (C\eps)^n.
\]
To prove the stronger bound stated in Theorem~\ref{thm:concentrationOff}, we do a ``second-moment'' type version of the above, which adds some extra complications, but traces the same shape as the above.

\section*{Acknowledgments.}
I would like to thank Marcus Michelen, for some useful discussion on topics related to the sphere packing literature. I would also like to thank Rob Morris and Marcus Michelen for comments on an earlier draft. Finally I would like to thank the volunteers working to put together the ICM surveys.

\bibliography{Bib.bib}
\bibliographystyle{siamplain}

\end{document}